\documentclass{article}

\makeatletter
\def\sigmabar{\overline{\sigma}}
\def\tauhat{\widehat{\tau}}
\def\bhat{\widehat{b}}
\def\chat{\widehat{c}}
\def\Zhat{\widehat{Z}}
\def\tautilde{\widetilde{\tau}}
\def\qq#1{\{#1\}}

\DeclareSymbolFont{AMSb}{U}{msb}{m}{n}
\DeclareMathSymbol{\N}{\mathbin}{AMSb}{"4E}
\DeclareMathSymbol{\Z}{\mathbin}{AMSb}{"5A}
\DeclareMathSymbol{\R}{\mathbin}{AMSb}{"52}
\DeclareMathSymbol{\Q}{\mathbin}{AMSb}{"51}
\DeclareMathSymbol{\I}{\mathbin}{AMSb}{"49}
\DeclareMathSymbol{\C}{\mathbin}{AMSb}{"43}

\newcommand{\arcsinh}{\mathop {\operator@font arcsinh}\nolimits}
\newcommand{\Cl}{\mathop {\operator@font Cl}\nolimits}
\makeatother

\begin{document}

\title{Some notes on the inverse problem for braids}

\author{Jonathan Fine\\Milton Keynes\\England\\\texttt{jfine@pytex.org}}

\date{13 August 2006}

\maketitle

\begin{abstract}
  The Kontsevich integral $Z$ associates to each braid $b$ (or more
  generally knot $k$) invariants $Z_i(b)$ lying in finite dimensional
  vector spaces, for $i = 0, 1, 2, \ldots$.  These values are not yet
  known, except in special cases.  The inverse problem is that of
  determining $b$ from its invariants $Z_i(b)$.
  
  In this paper we study the case of braids on two strands, which is
  already sufficient to produce interesting and unexpected
  mathematics.  In particular, we find connections with number theory,
  numerical analysis and field theory in physics.  However, we will
  carry this study out with an eye to the more general case of braids
  on $n$ strands.  We expect that solving the inverse problem even for
  $n=3$ will present real difficulties.  Most of the concepts in this
  paper also apply to knots, but to simplify the exposition we will
  rarely mention this.
  
  The organisation and bulk of the writing of this paper predates its
  most significant results.  We hope later to present better and
  develop further these results.

\end{abstract}

\section{Introduction}

The Kontsevich integral $Z$ associates to each knot $k$ in $\R^3$ (or
braid $b$ on $n$ strands) quantities $Z_i(k)$ for $i=0, 1, 2, \ldots$
lying in finite dimensional vector spaces $A_i$ (respectively $Z_i(b)$
lying in $A_{n,i}$).  These quantities depend only on the isotopy
class of the knot or braid.  This paper starts an investigation of the
following question: Suppose values $c_i$ in $A_i$ are given.  Is it
possible to find a weighted sum $k = \sum a_j k_j$ of knots $k_j$ such
that $c_i = \sum a_j Z_i(k_j)$ for all $i$?  In other words, can we
calculate the knot (or braid) from the values that result from
applying the Kontsevich integral.  We call this the \emph{inverse
  problem} (to that of calculating the Kontsevich integral).

In this paper we study the case of braids on two strands, which is
already sufficient to produce interesting and unexpected mathematics.
In particular, we find connections with number theory, numerical
analysis and field theory in physics.  However, we will carry this
study out with an eye to the more general case of braids on $n$
strands.  We expect that solving the inverse problem even for $n=3$
will present real difficulties.  Most of the concepts in this paper
also apply to knots, but to simplify the exposition we will rarely
mention this.

The subject matter of this paper is relatively elementary, and to make
it more accessible we have provided an informal exposition of those
parts of the Kontsevich integral that are required to motivate and
understand the calculations on $B_2$, which are the main results in
this paper. This will also help establish the point of view that leads
up to the calculations.  The full definition is required only for
braids on $3$ or more strands.

The interested reader should consult \cite{BN1} for a formal
exposition of the Kontsevich integral, \cite{BN2} and \cite{BN3} for a
problem and its solution that seem to be related to this paper, and
\cite{Ohtsuki} for a recent survey of open problems.  In addition,
\cite{Birman} provides a survey that encompasses both the old and new
points of view in knot theory.

A solution to the inverse problem will in general use infinite sums of
braids, and will therefore involve consideration of convergence.  For
most of this paper we take a naive approach, to provide examples that
might later inform a concluding discussion of this question.

The Kontsevich integral, considered algebraically, is somewhat novel.
Part of this paper is a study, motivated by the inverse problem, of
some of the properties of this new algebraic structure.

\subsection*{Note}

This is a preliminary version of this paper.  During its writing, the
relationship with the Dirichlet beta function was discovered.  We hope
in a future version to reorganise the material, and perhaps moving
some of the content to other articles.  In particular, we would like
to present earlier the calculations in sections~6 to~9, and move the
general definition of the Kontsevich integral to the end of the paper.

Comments are welcome.

\section{Braids on two strands}

The braid group $B_2$ on two strands is generated by the say clockwise
half-twist $\sigma_1$, which interchanges the first and second
strands.  Its inverse $\sigma_1^{-1}$ we will denote by $\sigmabar_1$.
For $B_2$, but not for higher braid groups, there are no relations
between the generators (other than $\sigma_1 \sigmabar_1$ being the
identity).  

Twists in $B_2$ can be \emph{counted}.  The Kontsevich integral
provides a means of \emph{measuring} the twist.  Each horizontal slice
through a braid has a certain amount of twist in it.  Now divide the
braid into many thin slices.  Summing the twist in all the slices, and
taking the limit to obtain the integral, gives $Z_1(b)$.  Doing the
same for the product of the twist in all unordered pairs of slices
gives $Z_2(b)$, and for unordered triples $Z_3(b)$.

By elementary calculus
\[
 1= 1, \int_1^1 1 dx = 1, \int_{0<x<y<1} 1 dx\,dy = 1/2, 
\int_{0<x<y<z<1} 1 dy\,dy\,dz = 1/6, \mbox{etc}
\]
and so the integrals that are the Kontsevich integral are readily
calculated when the twist is the same in each slice throught the
braid, as it is for $B_2$.  They are simply the volume of a simplex
(in a flat space).

For $B_2$ the value spaces $A_{2,i}$ are one-dimensional vector
spaces, and the difference $\tau_1 = \sigma_1 -\sigmabar_1$ defines a
non-zero element $t$ in $A_{2,1}$ (see~\S3). As an abstract
group, $B_2$ is isomorphic to $\Z$. In general, the group law induces
a non-commutative composition
\begin{equation}
  \label{A-comp}
  A_{d,i} \otimes A_{d,j} \to A_{d,i+j}
\end{equation}
and therefore $t$ induces elements, which we will denote by $t^i$, in
$A_{2,i}$.

With these identifications we have
\begin{equation}
 Z_i(\sigma_1) = \frac{1}{i!}\left(\frac{t}{2}\right)^i \in A_{2,i}
 \>,\qquad
 Z_i(\sigmabar_1)= \frac{1}{i!}\left(\frac{-t}{2}\right)^i \in A_{2,i}
\end{equation}
and so in this situation it is convenient to identify $\bigoplus
A_{2,i}$ with $\C[t]$ and then we can write
\begin{equation}
  Z(\sigma_1) = e^{t/2} \>, \qquad  Z(\sigmabar_1) = e^{-t/2}
\end{equation}
where the exponentials stand for power series.  That the product
$\sigma_1\sigmabar_1$ is the identity (or trivial) braid is then
represented by the equation
\[
    e^{t/2} \cdot e^{-t/2} = 1 \>.
\]

In this setting it can be helpful to write
\[
    q=\sigma_1\>,\qquad q^{-1}=\sigmabar_1
\]
and sometimes we will write $p$ for $q^{-1}$.  The Kontsevich integral
then consists of writing $q = e^{t/2}$, $p=e^{-t/2}$ and taking the
power series expansion in $t$.  Something similar can be done for the
Alexander and Jones polynomial invariants of knots.  They are usually
written as Laurent polynomials in $q$.  If $P_k(q)$ is such a
polynomial, then by a result of Lin the coefficient of $t^i$ in the
power series expansion of $P_k(e^t)$ (or $P_k(e^{t/2})$) is a linear
function of $Z_i(k)$.

For braids on two strands, the inverse problem consists of finding a
Laurent polynomial (in $q = e^{t/2}$) that has a given power series
expansion.  This is also part of the problem of computing, say, the
Alexander polynomial from the Kontsevich invariants.  The reader might
like now to look at \S6, where an inverse to the Kontsevich integral
for $B_2$ is calculated, and the following two sections, which
investigate its properties.  The immediately following sections
motivate these calculations.

\section{The value spaces $A_{n,i}$ and their properties}

The value spaces $A_{n,i}$ are not abstract vector spaces.  They
consist of formal sums of braids, modulo relations.  Here's how.  The
strands of a braid do not cross.  A \emph{braid with double points} is
like a braid, except the strands are allowed to cross transversally,
at double points.  Each double point can be resolved (to remove the
intersection) in two ways.  Thus, up to a sign, each braid with $j$
double points determines an alternating sum $b$ of $2^j$ (double point
free) braids.  Let $A_{n,(j)}$ denote the span of all such $b$
(arising from braids on $n$ strands, and with $j$ double points).
Vassiliev first developed this approach to knot invariants.

For example, $\tau_1 = \sigma_1 -\sigmabar_1$ arises from a braid with
a double point, and ${\tau_1}^2 = (\sigma_1 -\sigmabar_1)^2$ arises
from two double points.  Note that by `dropping a double point'
$A_{n,(i+1)}$ is a subset of $A_{n,(i)}$. For example,
\[
\tau_1^2 = \sigma_1(\sigma_1 -\sigmabar_1) -
\sigmabar_1(\sigma_1 -\sigmabar_1)
\]
is both an element of $A_{2,(2)}$ and the difference of two elements
in $A_{2,(1)}$.

The value spaces $A_{d,i}$ are the quotient spaces $A_{d,(i)} /
A_{d,(i+1)}$. This is how $\tau_1 = \sigma_1 - \sigmabar_1$ induces
the element $t$ in $A_{2,1}$.  The differences $2(\sigma_1 -e)$ and
$2(e-\sigmabar_1)$ also induce the same element $t$ (where here $e$
is the group identity).

When the sum of braids $b$ lies in $A_{d,(i)}$, the $i$-th order
Kontsevich invariant $Z_i(b)$ has the special property that it is
simply the residue of $b$ in $A_{d,i}$.  Because $A_{d,(i)} \subseteq
A_{d,(j+1)}$ for $j<i$, it follows from this special property that
$Z_j(b)$ is zero for $j<i$ when $b$ lies in $A_{n,i}$.

\section{Weak and strong inverses}

We can now return to the inverse problem.  Let us say that a weighted
sum $a = \sum a_i b_i$ of braids is \emph{focussed on order $r$} if
$Z_j(a) = \sum a_i Z_j(b_i)$ is zero for $j \neq r$.  A focussed sum
of braids is something like an eigenvector, with eigenvalue $r$.
Suppose that for each $c$ in the value space $A_{n,r}$ we can find a
sum $a$ of braids such that $a$ is focussed on order $r$ and also that
$Z_r(a)=c$.  This is a special case of the inverse problem and from
its solution, subject to convergence, we can solve the general
problem.  (If values $c_i$ are given, find $a_i$ such that $Z_j(a_i)$
is $c_i$ for $i=j$ and zero otherwise, and then set $a$ to be $\sum
a_i$.)

Suppose that $c$ in $A_{d,j}$ is given and that $a = \sum a_i b_i$ is
both focussed of order $j$ and a solution to $Z_j(a) = c$.  In other
words, $a$ solves the inverse problem for $c$.  From the residue
property it follows that $a$ lies first in $A_{d,(0)}$, then
$A_{d,(1)}$ and so on up to $A_{d,(j)}$.  Therefore, $a$ is
characterized by two properties.  First, it is a representative in
$A_{d,(j)}$ of the residue class $c$ in $A_{d,j}$.  Second, $Z_k(c)=0$
for all $k>j$.

We now come to one of the key definitions of this paper.  Suppose that
for each $i$ we are given a linear section $Y_i$ of the quotient map
$\pi_i:A_{d,(i)} \to A_{d,i}$. (This means that $Y_i$ selects a single
representative from each residue class.  This has some analogy with
Hodge theory, which selects a single representative in each de Rham
cohomology class.)

In that case we will say that $Y$ is a \emph{weak inverse} to the
Kontsevich integral~$Z$.  In terms of $Z$ it has the \emph{diagonal}
and \emph{lower triangular properties}
\begin{eqnarray}
Z_i(Y_i(Z_i(b))) &=& Z_i(b)  \label{diag-prop}\\
Z_j(Y_i(Z_i(b))) &=& 0 \qquad j < i \label{lower-t-prop}
\end{eqnarray}
while a \emph{strong inverse} has the additional \emph{upper
  triangular property}
\begin{eqnarray}
Z_j(Y_i(Z_i(b))) &=& 0 \qquad j > i \>. \label{upper-t-prop}
\end{eqnarray}

Although (\ref{diag-prop}) and (\ref{lower-t-prop}) are written in
terms of $Z$, that $Y$ is a weak inverse is a statement about $Y$ and
the value spaces $A_{n,i}$ alone.  However, the upper triangular
property (\ref{upper-t-prop}) that additional characterises a strong
inverse necessarily involves the Kontsevich integral $Z$.

For later use, we will say that $Y$ is \emph{coherent} if it respects
the product (\ref{A-comp}).  More exactly, this means that
\begin{equation}
  Y_{i+j}(c) = Y_{i}(a) Y_{j}(b)
\end{equation}
where $c$ is the product of elements $a$ and $b$ in the residue spaces
$A_{n,i}$ and $A_{n,j}$ induced by the group law in $B_n$, while the
product on the right side is a product of sums of elements of $B_n$.

\section{Strengthening a weak inverse}

A key problem in this area, which is still open, is whether the
Kontsevich invariants $Z_i$ distinguish braids (and knots).  Indeed,
this problem was a major reason for studying the inverse problem.  One
formulation is this: Define $A_{n,(\infty)}$ to be the intersection of
all the $A_{n,(i)}$.  Does $A_{n,(\infty)}$ contain any non-zero
elements?

If $A_{n,(\infty)}$ has non-zero elements, then there cannot be a
unique strong inverse $Y$, because any such can be modified to $Y +
\epsilon$, where each $\epsilon_i$ is a linear map from $A_{n,i}$ to
$A_{n,(\infty)}$.  

Conversely, any weak inverse $Y$ can be successively modified to
produce a sequence $Y=Y^{(0)}, Y^{(1)}, \ldots$ which, if convergent,
converges to a strong inverse~$Y^{(\infty)}$.  The next section gives
an example of this process.  Here we describe the process in general.
Suppose $Y=Y^{(m)}$ satisfies not only (\ref{diag-prop}) and
(\ref{lower-t-prop}) but also the special cases
\begin{equation}\label{r-vanish}
Z_j(Y_i(Z_i(b))) = 0 \qquad i < j < m
\end{equation}
of (\ref{upper-t-prop}).  Now define $Y'=Y^{(m+1)}$ by
\[
    Y'_i(v) = Y_i(v) - Y_m(Z_m(Y_i(v))
\]
for $i<m$.  In other words, $Y'$ is $Y$ corrected by the lift of the
amount that $Y$ is `off' at $j=m$.  As $Y_i'-Y_i$ takes values in
$A_{n,(m)}$, upon which $Z_i$ vanishes (because $i<m$), $Y'$ is still
a weak inverse satisfying~(\ref{r-vanish}).  However, writing $b'_i =
Y_i(Z_i(b))$,
\begin{eqnarray}
    Z_m(Y'_i(Z_i(b))) &=& Z_m(b'_i) - Z_m(Y_m(Z_m(b'_i)))\\
    &=&  Z_m(b'_i) - Z_m(b'_i) = 0
\end{eqnarray}
which shows that $Y'$ also satisfies (\ref{r-vanish}), but with $m$
replaced by $m+1$.

Provided the sequence $Y$, $Y''$, $\ldots$ converges, this process
produces from the weak inverse $Y$ a strong inverse~$Y^{(\infty)}$.
Moreover, provided $A_{d,(\infty)}$ has only the zero element, any
other so obtained sequence will either diverge, or converge to the
same limit.

However, it might not be easy to find a weak inverse that, when
strengthened as above, converges.  This is because computing the
Kontsevich integral is another key problem that is still open, and
from a strong inverse the Kontsevich integral can be computed.  (The
same is not true of a weak inverse.)

The key to solving this `inverse-inverse' problem is to exploit the
upper triangular property, together with the residue property for
elements of $A_{n,(i)}$.  Briefly, given a strong inverse $Y$ and a
braid $b$ one can use residues to compute $c_0=Z_0(b)$, and then
subtract $Y_0(c_0)$ from $b$ to obtain a sum of braids $b'$.  By
construction, $Z_i(b)$ and $Z_i(b')$ are equal for $i \geq 1$.  Now
repeat this process to calculate $Z_1(b)$, $Z_2(b)$ and so on.  This
is related to `actuality diagrams', as described in
Birman~\cite{Birman}.

\section{The $\tau$-seeded weak inverse for $B_2$ and its limit}

We will now compute a weak inverse for $B_2$, and also its sequence of
strengthened approximations.  We will compute a coherent inverse, and
so, as the value spaces $A_{2,i}$ for $B_2$ are generated in
degree~$1$, it is enough to choose a section of $\pi_1:A_{2,(1)} \to
A_{2,1}$, or in other words a sum of knots with double points, whose
residue is non-zero.

The three simplest elements of $B_2$, and their Kontsevich integrals, are
\begin{eqnarray}
1 &\mapsto& 1\\
\sigma_1 &\mapsto& 1 + 1/2 t + 1/8 t^2 + 1/48 t^3 + 1/384 t^4 \nonumber \\
                    &&\qquad + 1/3840 t^5 + 1/46080 t^6 + 1/645120 t^7\ldots\\
\sigmabar_1 &\mapsto& 1 - 1/2 t + 1/8 t^2 - 1/48 t^3 + 1/384 t^4 \nonumber\\
                    &&\qquad  - 1/3840 t^5 +  1/46080 t^6 - 1/645120 t^7\ldots
\end{eqnarray}
and in this section we choose to use $\tau_1 = \sigma_1 - \sigmabar_1$
as the weak inverse which, in this section, we will write as $\tau$.

To begin with we have
\begin{equation}
\tau \>\mapsto \> t + 1/24 t^3 + 1/1920 t^5 + 1/322560 t^7 \ldots
\end{equation}
which is `off' by $1/24t^3$. By coherence $t^3$ lifts to $\tau^3$ plus
higher order terms, and so $-1/24\tau^3$ is the correction.  The
strengthened approximation is then
\begin{equation}
\tau -1/24 \tau^3 \> \mapsto \> t +  (1/24-1/24)t^3 - 3/640 t^5 + \ldots
\end{equation}
where the correction not only cancels the $t^3$ but also changes the
coefficient of $t^5$ and higher order terms.  We are now `off' by
$-3/640 t^5$.  As before, by coherence $t^5$ lifts to $\tau^5$ plus
higher order terms, and so the next approximation is
\begin{equation}
    \tau -1/24 \tau^3 + 3/640 \tau^5
    \>\mapsto\> 
    t + (-3/640 + 3/640)t^5 + 5/7168 t^7 \ldots
\end{equation}
and by now the algorithm should be clear.

Here we have applied the general algorithm for strengthening a weak
inverse.  In the present case, the lifting is coherent, and generated
in degree one by a single generator (and so everything commutes).
This allows us to use classical inverse functions as a short-cut.

The mapping
\[
    t \>\mapsto\> e^{t/2} - e^{-t/2} = 2 \sinh\frac{t}{2}
\]
expresses $Z(\tau)$ as a function of $t$. Note that $Z$ is a linear
function of its argument, and so $Z(2\tau)$ is $4\sinh(t/2)$ and not
$2\sinh(t)$.  However $Z(\tau^2)$ is $(Z(\tau))^2$, because $Z$
respects the group law on $B_n$.  We will now show that the inverse
function
\[
    \tau \mapsto 2 \arcsinh \frac{\tau}{2}
\]
expresses the limit $Y^{(\infty)}(t)$ as a function of $\tau$.  

To do this, let $P_r(\tau)$ be the $r$-th order polynomial
approximation to $2\arcsinh(\tau/2)$, and then expand $P_r(2
\sinh(t/2))$ as a power series in $t$.  By the classical inverse
function theorem, this will be equal to $t$ plus terms of degree $>r$.
But $Z(P_r(\tau))$ is equal to $P_r(Z(\tau)) = P_r(2 \sinh(t/2))$.
Thus, $t \mapsto P_r(\tau)$ is an order~$r$ strengthening of our
original lifting.

The one-line Maple program
\begin{verbatim}
  series(2*arcsinh(tau/2), tau, 8)
\end{verbatim}
then gives us the following strong lifting
\begin{equation}
\label{tau-lift}
    t \mapsto \tau -\frac{1}{24}\tau^2 + \frac{3}{640} \tau^3
        -\frac{5}{7168} \tau^5 + \frac{35}{294912} \tau^7
        -\frac{63}{2883584} \tau^9 + \frac{231}{54525952} \tau^{11}
        - \ldots
\end{equation}
whose properties we will study in the next two sections.

\section{Some properties of the $\tau$-solution for $B_2$}

First, we note that (\ref{tau-lift}) produces a sum of powers of
$\tau$, and that $\tau$ is $\sigma_1 - \sigmabar_1$.  In the next
section we will expand (\ref{tau-lift}) as a sum of powers of
$\sigma_1$ and $\sigmabar_1$.  Here, we note that the rapidly
increasing denominators are likely to be helpful in establishing the
convergence of this sum.

For the remainder of this section, we look at some of the properties
of the numerators and denominators of (\ref{tau-lift}).  According to
Sloane's Online Encyclopedia of Integer Sequence (OEIC), the
numerators are sequence A055786, the numerators of the Taylor series
expansion of $\arcsinh(x)$.  I have not yet tried to prove this for
all~$n$.

The denominators are Sloane's sequence A002553, coefficients for
numerical differentiation.  Again, I have not yet tried to prove this
for all~$n$.  The Encyclopedia article, at the time of writing, does
not note this connection, and gives no formula for these coefficients.
We hope to  say more on this elsewhere.

These numbers also appear in theoretical physics.  One evening
Sloane's Encylopedia was down, and so I did a Google search instead.
To my surprise, I found not only cached copies of OEIC pages, but also
Czarnecki and Smirnov's paper \cite{CzarSmirn} on quantum field
theory.  Table~1 of that paper (page~7) contains a column of numbers
$c_n$ whose odd entries are precisely the coefficients of
(\ref{tau-lift}), multiplied by $\pi$ and without the alternating
sign.  The even numbers in the column $c_n$ are the coefficients of
the Taylor expansion of the closely related function $2\arccos(t/2)$.
Again, the authors were not aware of this connection.

Their paper \cite{CzarSmirn} is on ``the master two-loop propagator'',
and it is illustrated by a Feynmann diagram that is similar to the
diagram for $\tau_1$ in $B_2$.  It may be that the special physical
case they consider is precisely sufficient to allow the relevant part
of their problem to coincide with the problem we consider in this
paper.  Again, I have not tried to prove this coincidence of values,
for all $n$ (and I doubt that I have sufficient physics to attempt
this problem).

Table~1 in their paper also contains a column of numbers $a_n$.  These
have rather larger numerators and denominators than do the $c_n$, and
appear to be more complicated.  Czarnecki and Smirnov's formula (17)
and a result in their reference [12] taken together indicate that the
$a_n$ are part of an identity that involves $\zeta(3)$ and the Clausen
number $\Cl_2(\pi/3)$ (this is seen by putting $x=1$ in their (17)).
However, this may arise from physics that has nothing to do with
$B_2$.  Understanding the coincidence for their $c_n$ would, of
course, help greatly here.

\section{Expanding the $\tau$-inverse in terms of $\sigma_1$ and $\sigmabar_1$}

The successive truncations of the power series (\ref{tau-lift}) are
approximations to a strong inverse for the Kontsevich integral $Z$ on
the braid group $B_2$ on two strands, provided convergence is
satisfied.  In this section we begin to look at the question of
convergence.

First, we will express these approximations in terms of powers of
$\sigma_1$ and $\sigmabar_1$, which we will write as $q$ and $p$
respectively when it suits us.  The Maple procedure
\begin{verbatim}
    doit := proc(n)
    local tmp;
        tmp := convert(series(2*arcsinh(tau/2), tau, 2*n), polynom);
        expand(eval(tmp, tau=(q-1/q));
    end proc
\end{verbatim}
computes such approximations.

Here are the first few approximations, where we use $\qq{n}$ to denote
$(q^n-p^n)$.
\begin{eqnarray}
&&\qq{1}\nonumber\\[4pt]
&&\frac{9}{8}\qq{1} - \frac{1}{24}\qq{3}\nonumber\\[4pt]
&&\frac{1225}{1024}\qq{1} - \frac{245}{3072}\qq{3} 
    + \frac{49}{5120}\qq{5} - \frac{5}{7168}\qq{7} 
    \nonumber\\[4pt]
&&\frac{19845}{16384}\qq{1} - \frac{735}{8192}\qq{3}
    + \frac{567}{40960}\qq{5} - \frac{405}{229376}\qq{7} 
    + \frac{35}{294912}\qq{9}
    \nonumber\\
&&\frac{160083}{131072}\qq{1} - \frac{12705}{13107}\qq{3}
    + \frac{22869}{1310720}\qq{5} - \frac{5445}{1835008}\qq{7} 
    \nonumber\\
    &&\qquad\qquad\qquad\qquad
    + \frac{847}{2359296}\qq{9} - \frac{63}{2883584}\qq{11}
    \nonumber
\end{eqnarray}

The numerators of the coefficients of $\qq{1}$ are the denominators of
the odd terms in Wallis's approximation to $\pi/2$, and the
denominators are Wallis's numerators, times~$4$.  Thus, the
coefficient of $\qq{1}$ converges to $4/\pi$.

A Google search for 19845 and 16384 produces about ten relevant
references, most of which are concerned either with numerical analysis
(particularly wavelets) or with field physics.  For most of these
papers, many of the above coefficients of $\qq{n}$ also appear,
perhaps multiplied by a simple fraction.

A particularly useful paper found in this way is Fornber and Ghrist
\cite{FornGhrist}.  Its Table~5 in contains 6 rows, all of which agree
with the values given by the Maple formula.  They also give the
limiting values for the coefficients of $\qq{j}$ (assuming the values
here are always the same as theirs).  This value is $(-1)^{j+1}4/(\pi
j^2)$.

Assuming the numbers defined in this section continue to coincide with
those of \cite{FornGhrist} we thus obtain the following formula:
\begin{equation}
  \label{B2-inverse}
    \tauhat =
    \frac{4}{\pi}
    \left( 
        (q - p) 
        - \frac{1}{9}(q^3 - p^3) 
        + \frac{1}{25}(q^5 - p^5) 
        - \frac{1}{49}(q^7 - p^7)
        + \ldots
    \right)
\end{equation}
The author has done calculations that indicate that similar
expressions should exist for powers of $\tauhat$.

This indicates, but does not prove, that the strengthening process
applied to the weak inverse $\tau_1$ converges, to give
(\ref{B2-inverse}).

\section{Convergence and sums of braids}

Previously, we have adopted a naive approach to the question of
convergence.  Now that we have a candidate strong inverse
(\ref{B2-inverse}), it is time to be more careful.

To solve the inverse problem, we require three things of a sequence
$b_i$ of finite sums of braids:
\begin{enumerate}
\item That it converges to give a possibly infinite sum of braids.
\item That the Kontesevich integral of this infinite sum converges.
\item That the Kontsevich integral has desired value.
\end{enumerate}

We will calculate $Z_i(\tauhat)$.  First recall that
\[
    Z_i(q^j) =  \left(\frac{j}{2}\right)^i \frac{1}{i!},\qquad
    Z_i(p^j) = (-1)^i \left(\frac{j}{2}\right)^i \frac{1}{i!}
\]
and so $Z_{2i}(\tauhat)$ is zero.  This takes care of $0, 2, 4,
\ldots$.

Next note that
\[
   Z_1(\tauhat) = \frac{4}{\pi}\left(
     1 - \frac{1}{9}(3) + \frac{1}{25}(5) - \frac{1}{49}(7) 
     \right)
     = \frac{4}{\pi}(1 - 1/3 + 1/5 - 1/7 + \ldots )
\]
which is equal to $1$, as the last factor is Liebnitz's formula for
$\pi/4$.

We wish for $Z_i(\tauhat)$ to be zero, for $i$ equal to $3, 5, 7,
9, \ldots$.  However, this is equivalent to all of
\begin{eqnarray*}
  1^1 - 3^1 + 5^1 - 7^1 + 9^1 - \ldots \\
  1^3 - 3^3 + 5^3 - 7^3 + 9^3 - \ldots \\
  1^5 - 3^5 + 5^5 - 7^5 + 9^5 - \ldots \\
  1^7 - 3^7 + 5^7 - 7^7 + 9^7 - \ldots \\
  \mbox{and so on}
\end{eqnarray*}
being zero.  These sums diverge, but they are examples of the series
\[
\beta(s) = 1^{-s} - 3^{-s} + 5^{-s} - 7^{-s} + 9^{-s} - \ldots \\
\]
which is a formula for the Dirichlet beta function (which is an
example of a Dirichlet $L$-function).  In particular, the intriguing
equation
\[
2^{s-3} s! \pi Z_s(\tauhat) = \beta(2-s) \qquad \mbox{for $s=3,5,7,9,\ldots$}
\]
is suggested.

It is known that $\beta$ can, by analytic continuation, be extended to
the whole complex plane, and that $\beta(s)$ is zero for $s$ equal to
$-1, -3, -5, -7, \ldots$.  Therefore, $\tauhat$ is a strong inverse to
$Z$ for $B_2$, provided a satisfactory means is provided for dealing
with the divergent sums above.

\section{A balanced basis sequence inverse}

In the previous section we constructed a strong inverse for $B_2$,
subject to convergence, by taking a lifting $t \mapsto \tau$ of
$A_{2,1}$ to $_{2,(1)}$ and then using coherence to produce highter
order liftings.  Notice that in this process the lifting of $A_{2,1}$
to $A_{2,(1)}$ was successively changed.

In this section we will adopt a different approach.  We will produce a
sequence of braids $b_0, b_1, b_2, \ldots$ such that the residues of
$b_0, \ldots b_r$ form a basis for $A_{2,(0)} / A_{2,(r+1)}$.  Then,
given target values $c_i \in A_{2,i}$ we consider the linear system of
equations
\begin{equation}
    Z_i(a_0b_0 + \ldots + a_r b_r) = c_i \qquad i = 0, \ldots, r \>.
\end{equation}
Solving this equation for the $a$'s in terms of the $c$'s, and then
letting $r$ go to infinity, will produce a sequence of finite sums of
braids that will, subject to convergence, solve the inverse problem
for the target values $c$.

The putative strong inverse $\tauhat$ obtained in the previous section
can, of course, be obtained by making a suitable choice for the basis
sequence $b_0, b_1, b_2, \ldots$.  For example, the sequence $1,
\tauhat, \tauhat^2, \tauhat^3, \ldots$ would have that effect.  But
without knowing the answer ahead of time, it is not clear how to
define a basis sequence that is sure to produce the same inverse as
$\tauhat$.

In this section we use the basis sequence $1, \sigma, \sigmabar,
\sigma^2, \sigmabar^2, \ldots$ and we will do the calculations only
for the spaces $A_{1,0}, A_{1,2}, A_{1,4}, \ldots$. (In the next
section we will use $1, \sigma, \sigma^2, \ldots$)

In other words, we will be looking at the sequence of square matrices
\[
M_0 = \pmatrix{1},
M_1 = \pmatrix{1&1&1\cr0&1&-1\cr0&1&1},
M_2 = \pmatrix{1&1&1&1&1\cr0&1&-1&2&-2\cr0&1&1&4&4\cr0&1&-1&8&-8\cr0&1&1&16&16\cr}\\
\]
where the columns are the powers of $0, 1, -1, 2, -2, 3, -3, \ldots$,
together with their inverses $N_i = M_i^{-1}$.  (For simplicity, we
have dropped the $1/i!$.)

Calculations with Maple show us that the $(1,3)$ entries in $N_1, N_2,
\ldots$ are the sequence
\begin{equation}
    \label{zeta2-seq}
    -1, \frac{-5}{4}, \frac{-49}{36}, \frac{-205}{144}, 
    \frac{-5269}{3600},
    \frac{-5369}{3600},
    \frac{266681}{176400},
    \frac{-1077749}{705600},
    \ldots
\end{equation}
which are, up to a sign, the successive approximations to the infinite sum
\[
    1 + \frac{1}{4} + \frac{1}{25} + \frac{1}{36} + \frac{1}{49} + \ldots
\]
which converges to $\zeta(2) = \pi^2/6$.  The numerators and
denominators are in the OEIC.

Again using Maple, the $(1,5)$ entries of $N_2, N_3, \ldots$ are the sequence
\[
    \frac{1}{4},\>
    \frac{7}{18},\>
    \frac{91}{192},\>
    \frac{1529}{2880},\>
    \frac{37037}{64800},\>
    \frac{54613}{90720},\>
    \frac{63566689}{101606400},\>
    \ldots
\]
which has as differences between the terms
\[
    \frac{1}{4},\>
    \frac{5}{36},\>
    \frac{49}{576},\>
    \frac{41}{720},\>
    \frac{5269}{129600},\>
    \frac{767}{25200},\>
    \frac{266681}{11289600},\>
    \frac{1077749}{57153600},\>
    \ldots
\]
whose denominators (so far) divide the denominators in
(\ref{zeta2-seq}), and whose numerators are roughly the square of
those in (\ref{zeta2-seq}).

These calculations indicate that the $(1,3)$ terms are sums of the
reciprocals of squares, while the $(1,5)$ terms are sums of sums of
squares.

It seems from the evidence available that the inverse matrices $N_i$
will converge, and so subject to convergence will give an inverse
$\tautilde$ for $t$.  There is at present no evidence that $\tautilde$
and $\tauhat$ are equal.  This is clearly an important question.

(To be continued.)

\section{An unbalanced basis sequence inverse}

Here we look at the inverse due to the basis sequence
\[
    1, \sigma, \sigma^2, \sigma^3, \sigma^4, \ldots \>.
\]
As this sequence misses completely $\sigmabar, \sigmabar^2, \ldots$ it
is hardly possible that it can converge.

And indeed, using Maple, this is what we find.  However, some of the
sequences that arise are of interest, and appear in the OEIC.

(To be continued.)

\section{Convergence}

This section starts to discuss definitions that might allow one to say
that one has a strong inverse $Y$ to the Kontsevich integral $Z$, at
least in the case of $B_2$.  However, whether or not the definitions
work will also depend on the mathematical facts.

At present, we do not know if the limit $\tautilde$ of the inverses
due to the balanced basis sequence is equal to the $\tauhat$ inverse
defined via $2\arcsinh(x/2)$. 

For the rest of this section we will suppose that $b$ is a sequence
$b_1, b_2, b_3, \ldots$ of sums of braids, where each $b_i$ is a
\emph{finite} sum of braids.  Thus, there is no doubt that each $b_i$
is well defined.

We will say that the sequence $b$ is \emph{biconvergent} if
\begin{itemize}
\item[(a)] the $b_i$ converge to a possibly infinite sum of braids $\bhat$.
\item[(b)] the $Z_j(b_i)$ converge, for each $j$, to a value $\Zhat_j(b)$.
\end{itemize}
It may also be possible to define a concept of \emph{uniform}
biconvergence.

Condition (a) can also be written as $X_c(b_i)$ converges for all $c$,
where $c$ is an actual braid (not a sum of braids) and $X_c(b_i)$ is
the coefficient of $c$ in $b_i$.  Doing this treats $Z$ and $Y$ on
more of an equal footing.

In this paper we have shown (more or less) that $\tauhat$ is
biconvergent, and indicated that the limit $\tautilde$ of the balanced
basis inverses is also likely to be biconvergent.

Results along the following lines would allow the whole theory to hang
together.  Suppose $b$ and $c$ are biconvergent sequences of braids.
Then
\begin{enumerate}
\item If $\bhat = \chat$ then $\Zhat(b) = \Zhat(c)$.
\item The sequence $b+c$ is biconvergent.
\item $\widehat{(b+c)} = \bhat + \chat$.
\item $\Zhat(b+c) = \Zhat(b) + \Zhat(c)$.
\end{enumerate}

Properties (2)--(4) follow immediately from the definition of
biconvergence.  Property (1) has the following consequence. The
Kontsevich integral $Z$ cannot be directly evaluated on the infinite
sum $\tauhat$ defined in (\ref{B2-inverse}).  We saw that attempting
to do so would produce a divergent sum.  Property (1) asserts that
however the terms of (\ref{B2-inverse}) are rearranged, to ensure
(possibly uniform) convergence of the values $Z_j(\tauhat)$, the
limiting value is not changed.

The alternating harmonic sum $\sum(-1)^{i+1}/i$ multiplied by
$\sigma_1$ satisfies~(a) and (b) above, but its terms can be
rearranged to give any desired sum, as the series is not absolutely
convergent.  Therefore, we add to \emph{biconvergence} the condition
\begin{itemize}
\item[(c)] If $i<j$ then $b_i - b_j \in B_{n,(i)}$.
\end{itemize}

(To be continued.)

\section{Some comments on $B_3$}

For $B_2$, the local and global theories are basically the same,
because the twisting in each slice is constant.  For $B_3$, we night
be able to tackle the local problem without knowing how to evaluate
the Kontsevich integral.  This would be a study of three displaced
points in the plane.

Another approach to $B_3$ is to start with the Borromean subgroup --
remove any strand and the result is trivial.  What is this subgroup?
A strong inverse should not take up out of this group.  Therefore,
perhaps, a problem that can be solved.

\end{document}